\documentclass[10pt]{article}
\usepackage{amsmath}
\usepackage{amssymb}
\newtheorem{thm}{Theorem}[section]

\newtheorem{prop}[thm]{Proposition}

\newtheorem{prob}{Problem}
\newcommand{\g}{{\mathfrak g}}
\newcommand{\ghat}{{\hat{\mathfrak g}}}
\newcommand{\into}{\hookrightarrow}
\newcommand{\iso}{\simeq}
\renewcommand{\l}{\ell}
\newcommand{\La}{\left\langle}
\newcommand{\Ra}{\right\rangle}

\newcommand{\R}{{\mathbb R}}
\newcommand{\rank}{{\mathop{\mathrm{rank}}\nolimits}}
\newcommand{\pa}[1]{{\langle{#1}\rangle}}
\newcommand{\pf}{\noindent \textbf{Proof.} }
\newcommand{\qed}{\hfill $\blacksquare$ \medskip}
\newcommand{\tensor}{\otimes}
\newcommand{\w}{\omega}
\newcommand{\Z}{{\mathbb Z}}
\newcommand{\hdom}{{\begin{picture}(8,4)\multiput(0,0)(4,0){3}{\line(0,1){4}}%
\multiput(0,0)(0,4){2}{\line(1,0){8}}\end{picture}}}
\newcommand{\vdom}{{\begin{picture}(4,8)\multiput(0,0)(0,4){3}{\line(1,0){4}}%
\multiput(0,0)(4,0){2}{\line(0,1){8}}\end{picture}}}
\title{Embeddings of Schur functions into types $B/C/D$}
\author{Michael Kleber\thanks{Partially supported by an NSF Mathematical 
          Sciences Postdoctoral Research Fellowship.}\\
        Massachusetts Institute of Technology\\
        {\tt kleber@math.mit.edu}}
\date{5 October 2000 \\[1ex]
{\small Notice required by journal: {\em This work has been submitted
to Academic Press \\ for possible publication. Copyright may be
transferred without notice, after\\ which this version may no longer be
accessible.}}}
\begin{document}
\maketitle

\begin{abstract}
We consider the problem of embedding the semi-ring of Schur-positive
symmetric polynomials into its analogue for the classical types $B/C/D$.
If we preserve highest weights and add the additional Lie-theoretic
parity assumption that the weights in images of Schur functions lie in
a single translate of the root lattice, there are exactly two
solutions. These naturally extend the Kirillov--Reshetikhin
decompositions of representations of symplectic and orthogonal quantum
affine algebras $U_q(\ghat)$ (some still conjectural, some recently
proven).
\end{abstract}

\section{Introduction and Background}
\label{sec_intro}

Consider the infinite-dimensional vector space $Y$ over $\R$ whose
basis elements $\{v_\lambda\}$ are indexed by Young diagrams
$\lambda$, {\em i.e.} by all partitions of $n$ for all $n\geq0$.
There are two natural ring structures on $Y$, owing to the fact that
Young diagrams can be used to index the irreducible representations of
the various classical Lie groups.  One multiplication arises from the
decomposition of tensor products in type $A$; its structure constants
are the familiar Littlewood--Richardson coefficients.  The other
arises in the same way from tensor product decomposition in the other
classical types $B$, $C$ and $D$ (remarkably, all three series give
the same multiplication).

Our goal is to understand the embeddings of the former ring, $Y_A$,
into the latter, $Y_{BCD}$, in which
\begin{equation}
\label{map}
v_\lambda \mapsto v_\lambda + \sum_{\mu<\lambda} m_{\lambda\mu} v_\mu,
\qquad m_{\lambda\mu} \geq 0.
\end{equation}
Here $\mu<\lambda$ denotes the extended dominance order on partitions:
given $\lambda=\pa{\lambda_1,\dots,\lambda_r}$ and
$\mu=\pa{\mu_1,\dots,\mu_s}$, we say $\mu\leq\lambda$ if
$\mu_1+\cdots+\mu_k \leq \lambda_1+\cdots+\lambda_k$ for $1\leq k\leq
s$, where $\lambda_i=0$ for $i>r$.  (The order is ``extended'' because
we do not require that $|\lambda|=|\mu|$.)

The purpose of this paper is first to present a construction which
builds such an embedding out of any polynomial or formal power series
$p(x)\in\R[[x]]$ with the property that the symmetric function
$\kappa_p := \prod_i p(x_i) / \prod_{i<j} (1-x_i x_j)$ has a positive
expansion in terms of Schur functions.  The resulting embedding acts
on $v_\lambda$ by skewing $\lambda$ by $\kappa_p$, in a way we will
define precisely in Section~\ref{sec_thm}.  Second, we show that this
construction yields all of the desired embeddings.

Before we state our final and motivating result, we briefly clarify
the connection between $Y$ and representations of classical groups.
The Young diagrams $\lambda=\pa{\lambda_1,\ldots,\lambda_r}$ with
$\l(\lambda) := r \leq n$ rows index the dominant integral highest
weights, and therefore the irreducible finite-dimensional
representations, of $SL(n+1)$, $SO(2n+1)$, and $Sp(2n)$ (types $A$,
$B$ and $C$, respectively).  The situation in $SO(2n)$ (type $D$) is
slightly different; here they parametrize the restrictions to $SO(2n)$
of irreducible representations of $O(2n)$, which split into two
conjugate irreducibles if $\lambda_n>0$.

One remarkable property of this parameterization is that for fixed
Young diagrams $\lambda$, $\mu$, the multiset of diagrams giving the
summands in the tensor product of the $\lambda$ and $\mu$
representations is independent of $n$, provided $n$ is sufficiently
large (in particular, $n>\l(\lambda)+\l(\mu)$).  Moreover, in this
``stable limit'' the differences between types $B$, $C$ and $D$
vanish, giving us the ring structure we referred to above as
$Y_{BCD}$.  For a full discussion of these topics see the fundamental
work of Koike and Terada~\cite{KT}, to which we will refer in greater
detail later.

Returning to the type of embeddings described in equation~(\ref{map})
above, we can explain the constraints in the language of
representation theory.  Now our goal is to find a family of
representations $\{W(\lambda)\}$ of type $B/C/D$ which are a
homomorphic image of the irreducible representations $\{V(\lambda)\}$
of type $A$.  The embedding now describes the decomposition of
$W(\lambda)$ into irreducibles, and in the representation theoretic
setting we naturally demand that the constants $m_{\lambda\mu}$ are
nonnegative integers.  The Lie theoretic meaning of $\mu<\lambda$ is
that $\lambda-\mu$ is a sum of positive roots, guaranteeing that
$W(\lambda)$ has $\lambda$ as its maximal weight.  This agrees with
the above definition of the extended dominance order with the
additional constraint that $|\lambda| \equiv |\mu| \bmod 2$, to
require that the weights $\lambda$ and $\mu$ are in the same translate
of the root lattice.

We can now state our final result: given the Lie-theoretic integrality
and parity constraints, there are exactly two embeddings.  In the
above construction, they come from $p(x)=1$ and $p(x) = \frac1{1-x^2}
= 1+x^2+x^4+\cdots$.  The Schur function expansions of $\kappa_p$ is
as the sum of all $s_\lambda$, where $\lambda$ ranges over all
partitions which have only even length columns or rows, respectively.
The two embeddings are intertwined by the involution on $Y$ which
takes the basis vector $v_\lambda$ to $v_{\lambda'}$, where $\lambda'$ 
is the transpose Young diagram.

In the special case that $\lambda$ is a rectangle, the resulting
representations $W(\lambda)$ coincide with restrictions of certain
irreducible representations of the quantum affine algebra $U_q(\ghat)$
to $U_q(\g)$.  These decompositions were originally conjectured by
Kirillov and Reshetikhin~\cite{KR} and recently proved for
simply-laced classical $\g$ by Chari~\cite{Ch}; see
Section~\ref{sec_lie} for some details.  The two Lie-theoretic
embeddings give the decompositions for the orthogonal and symplectic
cases.  There is evidence to support the natural hope that our two
embeddings answer the same quantum groups question when $\lambda$ is
not rectangular.

The organization of the rest of this paper is as follows.  In
Section~\ref{sec_thm} we present the aforementioned construction, and
we state and discuss our main theorem (Theorem~\ref{thm_main}), which
asserts that this construction yields all of the embeddings we seek.
The proof of the main theorem is presented in Section~\ref{sec_pf}.
Finally, in Section~\ref{sec_lie} we discuss the connections to Lie
theory, and in particular the way in which these results generalize
the Kirillov--Reshetikhin representations of quantum affine algebras.

\medskip \noindent {\bf Acknowledgments. }
I am extremely grateful to Richard Stanley, Monica Vazirani, and
Vyjayanthi Chari for various conversations while this work was
underway, and to Nicolai Reshetikhin for my interest in the topic.
Thanks to Monica Vazirani as well for extensive comments on an early
version of this report.  The ACE computer package~\cite{ACE} was quite
helpful for many symmetric function computations.

\section{Constructing Embeddings}
\label{sec_thm}

In this section we work in the ring $\Lambda$ of symmetric functions
in countably many variables with rational coefficients (or, when
necessary, in a completion $\hat{\Lambda}$ where we allow infinite
formal sums of basis elements) and follow the standard notation of
Macdonald's book~\cite{Mac}.  Recall in particular that the Schur
functions $\{s_\lambda\}$, the ``universal characters'' of the
polynomial representations of $GL(n)$, are an orthonormal basis of
$\Lambda$ with respect to the standard inner product, and that the map
$\phi_A:Y\to\Lambda$ which takes $\lambda$ to $s_\lambda$ is an
embedding of the ring structure we called $Y_A$.  We use $\emptyset$
to denote the empty partition of zero; $s_\emptyset=1$ in $\Lambda$.

Extending the Schur function picture to the other classical types,
Koike and Terada~\cite{KT} defined two additional bases
$\{sp_\lambda\}$ and $\{o_\lambda\}$ of the same ring, consisting of
analogous universal characters for the symplectic and orthogonal
groups.  The maps $\phi_{C}:\lambda\mapsto sp_\lambda$ and
$\phi_{BD}:\lambda\mapsto o_\lambda$ are two different embeddings of
the ring structure we called $Y_{BCD}$.

We can use these embeddings to rephrase our quest for maps $Y_A \to
Y_{BCD}$ from equation~(\ref{map}) in terms of automorphisms of
$\Lambda$, with the minor inconvenience that we must break symmetry by
picking one of $\phi_C$ or $\phi_{BD}$.  We arbitrarily choose
$\phi_C$ in this and the next section.

\begin{thm}
\label{thm_constr}
Let $p(x)= 1 + a_1 x + a_2 x^2 + \cdots$ be a polynomial or formal
power series with $p(0)=1$, and define
$$
\kappa_p = \frac{\prod_{i=1}^\infty p(x_i)}
                {\prod_{1\leq i<j<\infty} (1-x_i x_j)} \in \hat{\Lambda}.
$$
Then the map $f_p:\Lambda\to\Lambda$ defined by
$$
f_p(s_\lambda) = \phi_C\circ\phi_A^{-1}(\kappa_p^\perp s_\lambda)
$$
is a ring homomorphism.  Here $\phi_C\circ\phi_A^{-1}$ is the map
$s_\mu \mapsto sp_\mu$, while $\kappa_p^\perp$ indicates skewing by
$\kappa_p$, the adjoint to multiplication by $\kappa_p$.
\end{thm}

\pf
We show that the map $f_p = \phi_C \circ \phi_A^{-1} \circ
\kappa_p^\perp$ is the automorphism of $\Lambda$ defined on the
homogeneous symmetric functions $h_n$ by $h_n \mapsto h_n + a_1
h_{n-1} + a_2 h_{n-2} +\cdots + a_{n-1} h_1 + a_n$.  This completely
defines the map, as the $h_n$ are an algebraically independent set of
generators for $\Lambda$.  Although $\kappa_p$ generally has an
infinite expansion in terms of Schur functions, recall that
$s_\mu^\perp s_\lambda=s_{\lambda/\mu}=0$ unless
$\mu\subseteq\lambda$, so $\kappa_p^\perp s_\lambda \in\Lambda$.
Thanks to R.~Stanley for showing us a special case of the following
argument.

We will rewrite $f_p$ using the inner product $\langle\cdot\,,\cdot
\rangle$ instead of $\phi_C \circ \phi_A^{-1}$ to transform Schur
functions to their symplectic analogues.  We will need two sets of
variables $(x_i)$ and $(y_j)$; the inner products are all taken with
respect to the $x$ variables.
\begin{align*}
f_p(s_\lambda(y))
 &= \sum_\mu \La \kappa_p^\perp s_\lambda(x)\,, 
       s_\mu(x) \Ra_{\!x}\, sp_\mu(y) \\
 &= \La s_\lambda(x)\,, \kappa_p 
       \left( \textstyle \sum_\mu s_\mu(x) \, sp_\mu(y) \right) \Ra_{\!x} \\
 &= \La s_\lambda(x)\,, \frac{\prod_i p(x_i)}{\prod_{i<j}(1-x_ix_j)} \,
       \frac{\prod_{i<j}(1-x_ix_j)}{\prod_{i,j}(1-x_iy_j)} \Ra_{\!\!x} \\
\intertext{%
The expansion $\sum_\mu s_\mu(x) \, sp_\mu(y)$ for
${\prod_{i<j}(1-x_ix_j)}/{\prod_{i,j}(1-x_iy_j)}$ is found
in~\cite{KT} based on formulas from Littlewood~\cite{Li}.  Cancelling,
we continue:
}
 &= \La s_\lambda(x)\,, \textstyle \prod_i p(x_i) \,
       \prod_i \left( \sum_{n\geq0} h_n(y) x_i^n \right) \Ra_{\!x} \\
 &= \La s_\lambda(x)\,, \textstyle \prod_i \sum_{n\geq0}
       (h_n(y) + a_1 h_{n-1}(y) + \cdots)\, x_i^n \Ra_{\!x} \\
 &= s_\lambda(y)\Big|_{h_n(y)\mapsto h_n(y)+a_1h_{n-1}(y)+\cdots}
\end{align*}
This last equality is the ring homomorphism $h_n(y)\mapsto
h_n(y)+a_1h_{n-1}(y)+\cdots$ applied to the Cauchy identity
$\prod_{i,j}(1-x_iy_j)^{-1}=\sum_\mu s_\mu(x)s_\mu(y)$.
\qed

\noindent{\bf Example \refstepcounter{thm}\thethm\label{ex}}
There are three noteworthy cases in which the Schur function expansion
of $\kappa_p$ was recognized by Littlewood (11.9 in~\cite{Li},
p. 238):
\begin{eqnarray*}
1/\prod_{i<j}(1-x_ix_j) &=& \sum_{\lambda} s_{(2\lambda)'} \\
1/\prod_i (1-x_i^2)\prod_{i<j}(1-x_ix_j) &=& \sum_{\lambda} s_{2\lambda} \\
1/\prod_i (1-x_i)\prod_{i<j}(1-x_ix_j) &=& \sum_{\lambda} s_{\lambda}
\end{eqnarray*}
That is, for $p(x)=1$, $p(x)=1/(1-x^2)$, and $p(x)=1/(1-x)$, we can
expand $\kappa_p$ as the sum of $s_\lambda$ for all $\lambda$ with
even column heights, all $\lambda$ with even row lengths, and all
$\lambda$, respectively.

When $p(x)=1$ and $f_p$ is the identity, this gives us the Character
Interrelation Theorem of Koike and Terada~\cite{KT}, which expands
$s_\lambda$ in the $sp$-basis.  For a concrete example, let us compute
$f_p$ on $s_\lambda$ for $\lambda=\pa{322}$.  Using the above
expansion, we see that $\kappa_p^\perp$ takes $s_\lambda$ to $\sum_\mu
s_{\lambda/\mu}$ where $\mu$ runs over all partitions with even column
heights.  Since $s_{\lambda/\mu}=0$ unless $\mu\subseteq\lambda$, the
only $\mu$ which contribute are $\pa{22}$, $\pa{11}$, and the empty
partition:
\begin{align*}
\kappa_p^\perp s_\pa{322} 
&= s_{\pa{322}/\pa{22}} + s_{\pa{322}/\pa{11}} + s_\pa{322} \\
&= (s_\pa{3}+s_\pa{21}) + (s_\pa{311}+s_\pa{221}) + s_\pa{322}
\end{align*}
So $s_\pa{322} = f_p(s_\pa{322}) = sp_\pa{322} + sp_\pa{311} +
sp_\pa{221} + sp_\pa{3} + sp_\pa{21}$.
\qed

We can now state our main theorem.

\smallskip\noindent{\bf Definition.}
Say $p(x)$ is {\em $\kappa$-positive} if $\kappa_p$ is $s$-positive,
{\em i.e.} has all nonnegative coefficients when written in the basis
of Schur functions. For instance, the above comments show that $1$,
$1/(1-x)$, and $1/(1-x^2)$ are all $\kappa$-positive.

\begin{thm}[Main Theorem]
\label{thm_main}
Let $f:\Lambda\to\Lambda$ be a ring homomorphism such that
$f(s_\lambda) = sp_\lambda + \sum_{\mu<\lambda} m_{\lambda\mu}
sp_\mu$ for some constants $m_{\lambda\mu}$.
\begin{enumerate}
\renewcommand{\labelenumi}{{\rm(\alph{enumi})}}
\item
If all $m_{\lambda\mu}\geq0$ then $f=f_p$ for some $\kappa$-positive
$p$, as in Theorem~\ref{thm_constr}.
\item
If we add the Lie-type assumption that all $m_{\lambda\mu}$ are
integers and $m_{\lambda\mu}=0$ unless $|\lambda|\equiv|\mu|\bmod2$,
then $p(x)$ is either $1$ or $1/(1-x^2)$.
\end{enumerate}
\end{thm}

We will defer the work of proving Theorem~\ref{thm_main} to the next
section.  For now we make the following straightforward
observations about part~(a):
\begin{enumerate}
\item
The converse is clear: if $p$ is $\kappa$-positive then $f_p$
certainly takes $s_\lambda$ to something $sp$-positive, since
$s_{\lambda/\mu}$ is always $s$-positive.
\item
If $f=f_p$ and all $m_{\lambda\mu}\geq0$ then $p$ is certainly
$\kappa$-positive: the only time $s_\emptyset$ appears in $s_\mu^\perp
s_\lambda$ is when $\lambda=\mu$, so the coefficient of $s_\lambda$ in
$\kappa_p$ is $m_{\lambda\emptyset}\geq0$.
\item
As a corollary of the theorem, all $m_{\lambda\mu}$ being nonnegative
implies $m_{\lambda\mu}$ is zero unless $\mu\subseteq\lambda$, a much
stronger condition than $\mu\leq\lambda$.
\end{enumerate}

Observations about part~(b) and connections to Lie theory will be made
in Section~\ref{sec_lie}.  In contrast with having only two solutions
in the Lie case, the full space of $\kappa$-positive polynomials and
formal power series is quite large and appears somewhat messy.

\begin{prob}
\label{prob_kpos}
Give a closed-form characterization of the $\kappa$-positive $p(x)$.
\end{prob}

Given $p(x)=1+a_1x+a_2x^2+\cdots$, the coefficient of $s_\lambda$ in
$\kappa_p$ is some polynomial in the $a_i$ for each $\lambda$, and it
seems {\em a priori} possible that no finite set of conditions is
equivalent to every one of those polynomials being positive.  This
problem, in other words, may have no good answer.  We make two
observations on $\kappa$-positivity.

\begin{prop}
\label{prop_spos}
$p(x)$ is $\kappa$-positive whenever $\prod_i p(x_i)$ is $s$-positive.
\end{prop}
This is clear because $\kappa_p$ is the above product multiplied by
the $s$-positive $\kappa_1$.  Stanley~\cite{GARSIA} observed that this
$s$-positivity condition is equivalent to a question of Schoenberg's
asking which $p$ are the generating functions for sequences giving
rise to totally positive infinite Toeplitz matrices.  The answer was
proved partially by Aissen, Schoenberg and Whitney~\cite{ASW} and
completed by Edrei~\cite{Edrei} and independently by
Thoma~\cite{Thoma}: $p$ must be of the form
$$
p(x) = e^{\gamma x} \textstyle 
  \prod_i (1+\alpha_i x) / \prod_j (1-\beta_j x)
$$
where $\alpha_i,\beta_j,\gamma$ are nonnegative real numbers and
$\sum_i\alpha_i$, $\sum_j\beta_j$ converge.  If $p(x)$ is a
polynomial, this is just demanding that all of its roots be real and
negative.  For an extended discussion see~\cite{EC2} (ex.~7.91,
pp.~481, 543ff).

The converse does not hold (consider $1/(1-x^2)$), and computational
evidence suggests that there are relatively open sets of $p(x)$ which
are $\kappa$-positive without meeting the $s$-positivity criterion.
Consider the generic quadratic $p(x)=1+bx+ax^2$.  For $a\geq a_0$ for
some $a_0$ (perhaps $a_0=1$), numerical evidence suggests that the
positivity of the coefficients of $s_\pa{2^t,1^t}$ in $\kappa_p$ for
all $t$ imply that $b^2\geq4a$, so $\kappa$- and $s$-positivity
coincide.  But for $a\leq a_1 \approx .39816\ldots$ (the real root of
$2z^3+3z^2+z-1$), the limiting coefficient seems to be that of
$s_\pa{32211}$, which only forces $b\geq\frac{a\sqrt{a+2}}{a+1}$.
\qed

\begin{prop}
\label{prop_dual}
The involution $p(x) \mapsto \frac1{(1-x^2)p(-x)}$ preserves
$\kappa$-positivity.
\end{prop}
Taking the transpose of Young diagrams commutes with both
multiplications ($Y_A$ and $Y_{BCD}$), and we will show it induces the
above involution on the set of $\kappa$-positive $p$.  Let $\omega$ be
the standard involution of $\Lambda$ (or $\hat{\Lambda}$) taking
$s_\lambda$ to $s_{\lambda'}$.  We will say that $p(x)$ and $q(x)$ are
{\em dual} if $\kappa_p=\omega\kappa_q$.  The expansions in terms of Schur
functions cited in Example~\ref{ex} show that $1$ and $1/(1-x^2)$ are
dual to one another, while $1/(1-x)$ is self-dual.

First note that $\omega$ acts on the symmetric functions $\prod_i
p(x_i)$ by $p(x)\mapsto1/p(-x)$; this can be derived after applying
$f_p$ to the well-known relation $(\sum h_n x^n)(\sum e_n (-x)^n)=1$.
As we observed above, $\kappa_p=\kappa_1\prod_i p(x_i)$, and we know
$\omega\kappa_1=\kappa_{1/(1-x^2)}$.  We conclude that the dual of
$p(x)$ is $1/(1-x^2)p(-x)$.

Note that we have shown that $\kappa_p^\perp$ takes $e_n$ to
$e_n+b_1e_{n-1}+b_2e_{n-2}+\cdots$, and thus that $f_p$ takes
$sp_\pa{1^n}$ to $sp_\pa{1^n}+b_1 sp_\pa{1^{n-1}}+b_2
sp_\pa{1^{n-2}}+\cdots$, where $q(x)=1+b_1x+b_2x^2+\cdots$ is 
dual to $p$.
\qed

Perhaps a classification of $\kappa$-positive $p(x)$ with integer
coefficients is possible.  These would describe automorphisms of
$\Lambda_\Z$, the ring of symmetric functions with integer
coefficients (which many authors just call $\Lambda$).  Given two
polynomials $r(x),t(x)\in\Z[x]$ with $r(0)=t(0)=1$ and all roots real
and negative, the $s$-positive $r(x)/t(-x)$ and its dual are
$\kappa$-positive; we do not know of any other integer examples.

\section{Proof of the Main Theorem}
\label{sec_pf}

In this section we present arguments and calculations to prove
Theorem~\ref{thm_main}.  The calculations involve finding relations
among various $m_{\lambda\mu}$, the coefficients of $sp_\mu$ in
$f(s_\lambda)$, generally by applying $f$ to the dual Jacobi--Trudi
formula expressing $s_\lambda$ as a determinant of a matrix of
elementary symmetric functions $e_n$.

These computations, naturally, will require expanding products $sp_\mu
sp_\nu$ into $\sum_\lambda d^\lambda_{\mu\nu} sp_\lambda$, where the
symplectic structure constants $d^\lambda_{\mu\nu}$ are not as
universally familiar as their type-$A$ analogues, the
Littlewood--Richardson coefficients $c^\lambda_{\mu\nu}$.  For what
follows it suffices to know that $d^\lambda_{\mu\nu} =
c^\lambda_{\mu\nu}$ when $|\lambda|=|\mu|+|\nu|$, and
$d^\lambda_{\mu\nu}=0$ unless $|\lambda|=|\mu|+|\nu|-2k$ for some
$k\in\Z_{\geq0}$.

The reader disinclined to descend into the depths of occasionally
cumbersome computation is invited to skip everything after the
statements of Propositions~\ref{prop_linear} and~\ref{prop_constant}
until their respective $\blacksquare$s.  The logic of the overall
argument should not be lost.

\begin{prop}
\label{prop_deg}
Let $f$ be a ring homomorphism such that
$f(s_\lambda)=sp_\lambda+\sum_{\mu<\lambda} m_{\lambda\mu}\,sp_\mu$.
We will denote $m_{\pa{1^i}\pa{1^j}}$ with $i>j$ by $m_{ij}$.

The $m_{ij}$ completely determine $f$, and each $m_{\lambda\mu}$ can be
expressed as a polynomial in the $m_{ij}$.  If we say $m_{ij}$ has
degree $i-j$, then each monomial of $m_{\lambda\mu}$ is of degree at
most $|\lambda|-|\mu|$ in the $m_{ij}$.
\end{prop}
When $\lambda=\pa{1^i}$ is a single column, the only $\mu<\lambda$ are
$\pa{1^j}$ for $j<i$, so knowing the $m_{ij}$ is equivalent to knowing
the image of $s_\pa{1^i}$ for all $i$.  (We use the dual form of the
Jacobi-Trudi identity so often precisely because there are so few
$\mu<\lambda$ when $\lambda$ is a single column, but not a single
row.)  Since $s_\pa{1^i}$ is the elementary symmetric function $e_i$
and the $e_i$ generate the ring $\Lambda$, the constants $m_{ij}$
determine the function $f$.

More concretely, the dual Jacobi--Trudi formula says
$s_\lambda=\det(e_{\lambda'_i-i+j})$, where
$\lambda'=\pa{\lambda'_1,\ldots,\lambda'_r}$ is the transpose of
$\lambda$, and the determinant is of an $r\times r$ matrix.  We apply
$f$ to both sides and find $m_{\lambda\mu}$ is the coefficient of
$sp_\mu$ in the determinant of the matrix whose $i,j$th entry is
$sp_\pa{1^k}+m_{k,k-1}sp_\pa{1^{k-1}}+\cdots+m_{k,0}sp_\emptyset$,
where $k=\lambda'_i-i+j$.  This is clearly a polynomial in the
$m_{ij}$.

If we assign to $m_{ij}$ degree $i-j$ and assign to $sp_\nu$ degree
$|\nu|$, then the $i,j$th matrix entry is homogeneous of degree
$\lambda'_i-i+j$.  This degree can only decrease when we expand
products, since $sp_\mu sp_\nu$ is a sum of $sp_\lambda$ with
$|\lambda|\leq|\mu|+|\nu|$.  So the terms of $f(s_\lambda)$ have
degree at most $|\lambda|$, and the coefficient $m_{\lambda\mu}$ is of
degree at most $|\lambda|-|\mu|$.
\qed

When we construct $f_p$ from a formal power series $p(x)$ as in
Theorem~\ref{thm_constr}, we get $m_{ij}=b_{i-j}$, where
$q(x)=1+b_1x+b_2x^2+\cdots$ is the dual of $p(x)$ in the sense of
Proposition~\ref{prop_dual}.  Therefore Theorem~\ref{thm_main}(a) is
precisely the claim that $m_{ij}$ depends only on $i-j$.

We proceed by induction on $d=i-j$.  Fix some integer $d>0$ for the
remainder of this proof.  Take some $f$ determined by its $m_{ij}$,
and suppose $m_{ij}$ is some constant $b_{i-j}$ whenever $i-j<d$.  Our 
task is to show that $m_{k,k-d}$ is independent of $k$.

To this end, construct $f_p$ from $p(x)$, the dual of the finite
polynomial $q(x)=1+b_1x+\cdots+b_{d-1}x^{d-1}$, as above.  This $f_p$
will also remain fixed for the duration; we can think of it as the
degree-$d$ polynomial truncation of our arbitrary embedding $f$.  Let
$m_{\lambda\mu}^{(p)}$ denote the constants describing $f_p$, just as
$m_{\lambda\mu}$ describe $f$.  By Proposition~\ref{prop_deg}, we
know that $m_{\lambda\mu}^{(p)}=m_{\lambda\mu}$ whenever
$|\lambda|-|\mu|<d$.

\begin{prop}
\label{prop_linear}
We compute some values of $m_{\lambda\mu}$: for $d\geq1$ and $k\geq d+2$,
\begin{enumerate}
\item
$m_{\pa{k},\pa{1^{k-d}}}=(-1)^{k-1}(m_{k,k-d}-2m_{k-1,k-1-d}+m_{k-2,k-2-d})$.
\item
$m_{\pa{k-1,k-1},\pa{k-2,1^{k-d}}} = -m_{\pa{k},\pa{1^{k-d}}}$.
\end{enumerate}
Since all $m_{\lambda\mu}$ are required to be nonnegative, we conclude 
that $m_{k,k-d}$ for fixed $d$ is a linear function of $k$.
\end{prop}
We begin with part~1, which is most of the work.  First note that we
want to compute $m_{\lambda\mu}$ where $\mu\nsubseteq\lambda$, so
$m_{\lambda\mu}^{(p)}=0$.  Also, since $|\lambda|-|\mu|=d$, we learn
from Proposition~\ref{prop_deg} that $m_{\lambda\mu}$ has degree at
most $d$, but any term of $m_{\lambda\mu}-m_{\lambda\mu}^{(p)}$ must
mention some $m_{ij}$ where $i-j\geq d$.  Combining these, we find
that $m_{\lambda\mu}$ must be a linear combination of $m_{i,i-d}$ for
various values of $i$.

%
%

Now consider the $k\times k$ matrix whose $i,j$th entry is
$f(s_\pa{1^{1-i+j}})$, whose determinant is $f(s_\pa{k})$ according to
dual Jacobi--Trudi.  We need only keep track of determinant
contributions linear in the $m_{i,i-d}$, all of which arise from
picking a permutation and from each matrix entry choosing either the
top-degree summand $sp_\pa{1^{1-i+j}}$ ($k-1$ times) or the
lower-degree summand $m_{1-i+j,1-i+j-d}sp_\pa{1^{1-i+j-d}}$ (once).
In this matrix, we have ones immediately below the main diagonal, and
zeros everywhere below those.  Thus the only permutations contributing
to the determinant correspond to choosing a composition
$k_1+k_2+\cdots+k_r=k$, as follows: within square diagonal blocks of
successive sizes $k_i$, choose all the subdiagonal unit entries and
the entry in the upper-right corner of the block.  The sign of this
permutation is $\prod(-1)^{k_i-1}=(-1)^{k-r}$.  From each
such term, our desired $m_{\lambda\mu}$ picks up a contribution of
$\pm m_{k_i,k_i-d}$ from the block of size $k_i$, for each $1\leq
i\leq r$.

Explaining the proposition's coefficients of $\pm1,\mp2,\pm1$ for a
distinguished block of size $k, k-1, k-2$ respectively is
straightforward.  The $m_{k,k-d}$ comes from the unique composition
with one part of size $k$, and gets sign $(-1)^{k-1}$, while the
$2m_{k-1,k-1-d}$ come from compositions $(k-1)+1$ and $1+(k-1)$ and
have sign $(-1)^k{-2}$.  When our distinguished block has size $k-2$
there are five compositions, which we can abbreviate as $\star11$,
$1{\star}1$, $11\star$, $\star2$ and $2\star$, where the $\star$ is
the block providing the $m_{k-2,k-2-d}$.  In three cases this has sign
$(-1)^{k-3}$ and in two cases $(-1)^{k-2}$, and the desired
coefficient is obtained.

This computation for $k-2$ makes it clear that each composition of $s$
with $t$ parts will contribute $t+1$ terms $m_{k-s,k-s-d}$,
corresponding to $t+1$ ways to insert a $\star$ designating a block of
size $k-s$, with sign depending on the parity of $t$.  To show that
these all cancel, it suffices to check that $\sum_t (-1)^t (t+1)
\mathop{\mathrm{comp}}(s,t)=0$ for any $s\geq3$.  Here
$\mathop{\mathrm{comp}}(s,t)$ denotes the number of compositions of
$s$ with exactly $t$ parts, equal to $\binom{s-1}{t-1}$, and the
resulting identity on binomial coefficients is easily verified by
induction.

Part~2 follows from part~1 by observing that $s_\pa{k-1,k-1}$ can be
computed via the (normal) Jacobi--Trudi formula as
$s_\pa{k-1}s_\pa{k-1} - s_\pa{k}s_\pa{k-2}$.  The same argument as
above shows that the desired coefficient
$m_{\pa{k-1,k-1},\pa{k-2,1^{k-d}}}$ is a linear combination of the
$m_{i,i-d}$.  These can arise in $f(s_\pa{k-1})f(s_\pa{k-1}) -
f(s_\pa{k})f(s_\pa{k-2})$ only by multiplying a coefficient linear in the
$m_{i,i-d}$ from one factor by a coefficient independent of the
$m_{i,i-d}$ and of maximal degree from the other.  On coefficients
independent of the $m_{i,i-d}$, $f$ agrees with $f_p$, so in
particular $m_{\lambda\mu}=0$ unless $\mu\subseteq\lambda$, and the
only term of maximal degree is $sp_\lambda$ itself.  Finally, we
observe that $s_\pa{k-2}$ is the only $s_\lambda$ mentioned here which
fits inside $\pa{k-2,1^{k-d}}$, thus we can only find the desired
shape in (minus) the product of $sp_\pa{k-2}$ with
$m_{\pa{k},\pa{1^{k-d}}} sp_\pa{1^{k-d}}$.

Thus $m_{\pa{k-1,k-1},\pa{k-2,1^{k-d}}}$ is $-m_{\pa{k},\pa{1^{k-d}}}$,
the negative of the coefficient calculated in part~1.  We conclude
that $m_{k,k-d}-2m_{k-1,k-1-d}+m_{k-2,k-2-d}$ must be zero for all
$k\geq d+2$, and therefore that $m_{k,k-d}$ depends linearly on $k$, as
claimed.  \qed

Now let us write the linear function $m_{j+d,j}$ as $\alpha_d +
j\beta_d$ for some constants $\alpha_d$ and $\beta_d$; we will
complete our induction by showing $\beta_d=0$.  Certainly
$\beta_d\geq0$, since $m_{j+d,j}$ must be nonnegative for arbitrarily
large values of $j$.  The following is essentially the same argument
with a little extra computation required.

\begin{prop}
\label{prop_constant}
We compute that $m_{\pa{k^d},\pa{(k-1)^d}} = \alpha_d - (k-1)
\beta_d$.  Since this must be nonnegative for arbitrarily large values
of $k$, we conclude $\beta_d=0$.
\end{prop}
Using the same strategy as before, we apply $f$ to the dual
Jacobi--Trudi matrix for $s_\pa{k^d}$.  It is again the case that
$m_{\pa{k^d},\pa{(k-1)^d}}^{(p)}$ is zero: while $\pa{(k-1)^d}$ is
indeed contained in $\pa{k^d}$, it only arises as $s_\mu^\perp
s_\pa{k^d}$ for $\mu=\pa{1^d}$, and we chose $p$ dual to a degree
$d-1$ polynomial, ensuring that $\langle s_\pa{1^d},\kappa_p\rangle=0$.
Thus the same degree argument as before shows that the desired
coefficient is a linear combination of the $m_{j+d,j}$.

Once again, we know that each contribution comes from taking some term
in the determinant of the Jacobi--Trudi matrix and replacing a single
factor $e_r$ with $m_{r,r-d}\,e_{r-d}$.  Note that $e_r$ must be
replaced in this way if it is above the main diagonal: any
superdiagonal $e_r$ not so replaced has $r>d$ and $s_\pa{(k-1)^d}$
cannot possibly appear.  Let us consider whether this replacement can
take place at the $i,j$ position in the matrix, where we would replace
$e_{d+j-i}$ with some multiple of $e_{j-i}$.  Certainly we need $j\geq
i$.

Now observe that the cofactor of the $i,j$ position in the matrix is
precisely the dual Jacobi--Trudi matrix for the skew shape
$\pa{(k-1)^d,i-1}/\pa{j-1}$.  The only $s_\lambda$ which appears
in the expansion of this skew shape and which is contained in
$\pa{(k-1)^d}$ is $\lambda=\pa{(k-1)^{d-1},k-1-j+i}$ --- that is,
$\pa{(k-1)^d}$ with a horizontal strip of length $j-i$ removed from
its last row.  We must multiply this by $e_{j-i}$ to get a factor of
$\pa{(k-1)^d}$, but the dual Pieri rule says multiplication by
$e_{j-i}$ adds a vertical strip of length $j-i$.  The only way this
vertical strip can fill the horizontal hole is if its length $j-i$ is
zero or one.

When $j-i=0$, we are looking at contributions from the identity
permutation: there are $k$ places to replace $e_d$ with $m_{d,0}e_0$.
When $j-i=1$ our permutation is one of the $k-1$ adjacent
transpositions, and we replace $e_{d+1}$ with $m_{d+1,1}e_1$.  Thus
$m_{\pa{k^d},\pa{(k-1)^d}} = k\,m_{d,0}-(k-1)m_{d+1,1}$.  Writing
$m_{j+d,j}$ as $\alpha_d + j\beta_d$, the coefficient is therefore
$\alpha_d-(k-1)\beta_d$.
\qed

Propositions~\ref{prop_linear} and~\ref{prop_constant} show that
$m_{ij}$ is a function only of $i-j$.  All such embeddings $f$ are
constructed in Theorem~\ref{thm_constr}, and the proof of
Theorem~\ref{thm_main}(a) is complete.

Theorem~\ref{thm_main}~(b) is now quite simple.  The Lie-theoretic
parity assumption means $p(x)$ is an even power series
$1+a_2x^2+a_4x^4+\cdots$.  We quickly compute that
$m_{\pa{2k+1,1},\emptyset}=a_{2k}-a_{2k+2}$ for all $k\geq0$, which
follows at once from $sp_\pa{i} sp_\pa{j}$ containing $sp_\emptyset$
once if $i=j$, and not at all otherwise.  Thus we have $1\geq a_2\geq
a_4\geq\cdots\geq0$.  The Lie-theoretic assumption that $a_2$ is an
integer leaves us with two cases.  If $a_2=0$ then all $a_i=0$ and
$p(x)=1$.  If $a_2=1$ then its dual $q(x)=1/(1-x^2)p(-x)$ is
$1+O(x^4)$.  But $q(x)$ is also an even function with integral power
series, so by the above, $q(x)=1$ and $p(x)=1/(1-x^2)$.  There is a
unique dual pair of solutions, as claimed.

\section{Connections to Lie Theory}
\label{sec_lie}

The motivation for this paper was to generalize some work of Kirillov
and Reshetikhin~\cite{KR} on the characters of certain irreducible
finite-dimensional representations of quantum affine algebras.
According to Theorem~\ref{thm_main}(b), there is a unique dual pair of
maps $Y_A \to Y_{BCD}$ with certain properties, which we mentioned
have a Lie-theoretic origin.  In this section we will explain how
those two maps give a generalization of the Kirillov--Reshetikhin
characters for the symplectic and orthogonal quantum affine algebras.

The paper in question concerned representation of the quantum affine
algebra $U_q(\ghat)$, where $\g$ is a simple Lie algebra of classical
type, and $\ghat$ is its corresponding affine Lie algebra.
Finite-dimensional representations of $U_q(\ghat)$ are not yet
well-understood.  Due to the embedding $U_q(\g)\into U_q(\ghat)$,
finite-dimensional $U_q(\ghat)$ modules can be said to have weights,
and any representation of $U_q(\ghat)$ can be viewed as a module over
the semisimple $U_q(\g)$, and thus decomposed into a direct sum of
$U_q(\g)$-irreducibles.  When $\g$ is of type $A_n$, irreducibles of
$U_q(\ghat)$ remain irreducible after restriction, but for types
$B/C/D$ this decomposition of irreducibles is nontrivial and not
known in general.

An earlier work of Kirillov~\cite{Ki} dealt with $\g$ of type $A$ and
investigated the decomposition of tensor products of ``rectangular''
irreducible representations --- that is, representations whose highest
weight is a multiple of a fundamental weight, whose Young diagram is a
rectangle.  The Bethe Ansatz and associated methods from mathematical
physics led to a formula for the decomposition; a central step in the
proof was the representation-theoretic identity
\begin{equation}
\label{eq_quad}
V(m\w_\l)^{\tensor2} \iso
 V((m+1)\,\w_\l) \tensor V((m-1)\,\w_\l) \oplus
 V(m\w_{\l-1}) \tensor V(m\w_{\l+1}).
\end{equation}
Here $V(m\w_\l)$ is the irreducible representation of $\g$ with
highest weight $m\w_\l$, where $m\in\Z_{\geq0}$ and
$\w_\l$ for $1\leq\l\leq\rank(\g)$ are the fundamental weights of
$\g$ (so $m\w_\l$ corresponds to a rectangle of $m$ columns each
of height $\l$).

The work of Kirillov and Reshetikhin~\cite{KR} applied the same
approach to the other classical Lie algebras, where the above identity
is no longer true for irreducible representations.  Instead it holds
when we replace the irreducible $V(m\w_\l)$ with certain reducible
representation $W(m\w_\l)$, which, according to the mathematical
physics origins, ought to be the decomposition of an irreducible
$U_q(\ghat)$ module, as discussed above.  Rephrasing their formulas in
the language of Young diagrams, when $\g$ is of type $B_n$ and $\l\leq
n-1$ or $D_n$ and $\l\leq n-2$, $W(m\w_\l)=\sum_\lambda
V(\lambda)$ summing over all shapes $\lambda$ which can be obtained
from the $\l\times m$ rectangle by removing vertical $2\times1$
dominos; for type $C_n$ and $\l\leq n-1$ we do the same but remove
horizontal dominos instead.

The identity~(\ref{eq_quad}) holds when $\w_{\l-1}$ and $\w_{\l+1}$
denote the two weights with the same length as and adjacent to $\w_\l$
in the Dynkin diagram of $\g$; the paper~\cite{KR} also gave a version
for when the $\w_\l$ node is trivalent or is adjacent to a shorter or
longer root, and corresponding values of $W(m\w_\l)$ for the $\l$ not
mentioned above, which we omit here.  The connection to $U_q(\ghat)$
was conjectural, as it is based on some unproven properties widely
believed to hold for the Bethe Ansatz.  Vyjayanthi Chari has recently
announced a proof~\cite{Ch} using entirely different techniques, in
the case that $\g$ is simply-laced and the weight $\w_\l$ appears with
multiplicity at most 2 in the highest weight of $\g$.  This proves the
Lie theory connection when $\g$ is of type $A$ and $D$, and for for
about half the choices of $\w_\l$ for exceptional $\g$ of type $E$.

The result is quite remarkable: the identity~(\ref{eq_quad}) implies
that the $W(m\w_\l)$ are given in terms of the fundamental
$W(\w_\l)$ by the dual Jacobi-Trudi formula, even when $\g$ is not
of type $A$.  A result of the author~\cite{K:poly} shows that this
requirement is so strong that there is in fact a {\em unique} family
of representations $W(m\w_\l)$ which obey the identity (including
the trivalent and different-length root extensions).

This brings us to the connection with the present work.  For $\g$ of
type $B/C/D$, the Kirillov--Reshetikhin representations $W(m\w_\l)$
(for $\l\ll\rank(\g)$) are indeed the restrictions of irreducible
$U_q(\ghat)$ representations if, and only if, those quantum affine
irreducibles obey the type-$A$ Jacobi-Trudi formula.  We extend this
proposition beyond the case of rectangles.

Theorem~\ref{thm_main}(b) gave us two distinguished automorphisms,
which we will now call $f_\hdom$ and $f_\vdom$, on the vector space
$Y$ whose basis vectors are labelled by Young diagrams.  These names
are mnemonic: $f_\hdom$ arises from skewing by $\sum_\lambda
s_\lambda$ where $\lambda$ ranges over all partitions will all even
rows, {\em i.e.} partitions built of horizontal dominos, and $f_\vdom$
similarly for vertical dominos.  In the notation of
Sections~\ref{sec_thm} and~\ref{sec_pf}, these are pulled back from
$f_p$ for $p(x)=1/(1-x^2)$ and $p(x)=1$, respectively, but this
asymmetry is misleading: if we had arbitrarily chosen
$\phi_{BD}:\lambda\mapsto o_\lambda$ instead of $\phi_C:\lambda\mapsto
sp_\lambda$ to embed our computations in the ring $\Lambda$, the
two $p(x)$ would be reversed.

\begin{prop}
Let $\g$ be of type $B/C/D$ and take some $\lambda$ with fewer than
$\rank(\g)-2$ parts.  Define a representation $W(\lambda)$ of $\g$ as
follows:
\begin{itemize}
\item
For symplectic $\g$ (type $C$), let $W(\lambda)$ have character
$\phi_C \circ f_\hdom(\lambda)$. \par
Equivalently, apply $o_\mu\mapsto sp_\mu$ to $s_\lambda$.
\item
For orthogonal $\g$ (type $B/D$), let $W(\lambda)$ have character
$\phi_{BD} \circ f_\vdom(\lambda)$. \par
Equivalently, apply $sp_\mu\mapsto o_\mu$ to $s_\lambda$.
\end{itemize}
These $W(\lambda)$ agree with the Kirillov--Reshetikhin values when
$\lambda$ is a multiple of a fundamental weight.  If irreducible
representations of quantum affine algebras obey type-$A$ algebraic
relations, then $W(\lambda)$ is isomorphic to an irreducible
representation of $U_q(\ghat)$ viewed as a $U_q(\g)$-module.
\end{prop}
Certainly an irreducible representation of $U_q(\ghat)$ viewed as a
$U_q(\g)$-module must have the properties we demanded in
Theorem~\ref{thm_main}: the $m_{\lambda\mu}$ are multiplicities of
$U_q(\g)$-irreducibles, so must be nonnegative integers, and the
weights $\mu$ that appear must have $\mu\leq\lambda$ in the
Lie-theoretic sense, from considering the lowering operators in
$U_q(\ghat)$.  The equivalent formulas for the character follow
because $1/(1-x^2)$ is dual to the identity on $\Lambda$ and
$\omega(sp_\mu)=o_{\mu'}$.  Note that $o_\mu\mapsto sp_\mu$ is the
ring homomorphism $h_i \mapsto h_i+h_{i-2}+h_{i-4}+\cdots$, and
$sp_\mu\mapsto o_\mu$ is its conjugate $e_i \mapsto
e_i+e_{i-2}+e_{i-4}+\cdots$.

To compare the Kirillov--Reshetikhin values, recall that if $\lambda$
is a rectangle, then $s_{\lambda/\mu}=s_\nu$ where $\nu$ is the skew
shape $\lambda/\mu$ rotated $180^\circ$.  \qed

For instance, our calculations in Example~\ref{ex} translated into the 
language of Lie theory would correspond to a type $B/D$ decomposition
$W(\w_1+2\w_3) \simeq V(\w_1+2\w_3) \oplus V(2\w_1+\w_3) \oplus
V(\w_2+\w_3) \oplus V(3\w_1) \oplus V(\w_1+\w_2)$.

\begin{prob}
Prove that these $W(\lambda)$ are the $U_q(\g)$-restrictions of
irreducible representations of $U_q(\ghat)$.
\end{prob}

Note that since we work in the stable limit where $n$ is assumed to be
sufficiently large that it is irrelevant, our results do not propose
decompositions of $W(\lambda)$ when $\lambda$ is supported on the spin
weights.

Until the recent announcement of~\cite{Ch}, there were almost no known
decompositions to compare with our proposed values of $W(\lambda)$.
The methods used there offer a way to calculate upper bounds on the
multiplicities in the decomposition of a certain canonical irreducible
representation of $U_q(\ghat)$ (the ``minimal affinization'' of Chari
and Pressley; see~\cite{ChP}) for simply-laced $\g$.  Vyjayanthi Chari
has kindly verified (private communication) that these upper bounds
coincide exactly with our $W(\lambda)$ for some special cases of type
$D$, {\em e.g.}  $\lambda=a\w_1+b\w_3$ and $\lambda=\w_2+\w_4$ (where
the decomposition is not multiplicity-free).

\end{document}